\documentclass[12pt, letterpaper]{article}
\usepackage[margin = 1in]{geometry}
\usepackage{graphicx,color}
\usepackage[T1]{fontenc}
\usepackage[english]{babel}
\usepackage[f]{esvect}
\usepackage{amsmath, amsthm, amssymb}
\usepackage{rotating}
\usepackage{mathrsfs}
\usepackage{url}
\usepackage{cite}

\def\tablenotes{\bgroup\parfillskip=0pt plus 1fil
\leftskip=0pt\relax \rightskip=0pt
\vskip2pt\footnotesize}
\def\endtablenotes{\vskip1pt\egroup}

\newtheorem{theorem}{Theorem}[section]

\newtheorem{definition}[theorem]{Definition}

\newcommand{\sinc}{\mathrm{sinc}}

\allowdisplaybreaks

\begin{document}

\title{Computation of Energy Eigenvalues of the Anharmonic Coulombic Potential with Irregular Singularities}

\author{M. Essaouini$^{\dag}$, B. Abouzaid$^{\ddag}$, P. Gaudreau$^{\P}$ and H. Safouhi$^{\S}$\footnote{Corresponding author. \newline The corresponding author acknowledges the financial support for this research by the Natural Sciences and Engineering Research Council of Canada~(NSERC) - Grant RGPIN-2016-04317.}\\
\\
$^{\dag}${\it Department of Mathematics, Faculty of Sciences}\\
{\it Chouaib Doukkali University, El Jadida, MA}\\
{\it essaouini.m@ucd.ac.ma}\\
$^{\ddag}${\it \'Ecole Nationale des Sciences Appliqu\'ees}\\
{\it Laboratoire des sciences de l'ing\'enieur pour l'\'energie}\\
{\it El Jadida, MA 24002}\\
{\it bhabouzaid@yahoo.fr}\\
$^{\P}${\it Facult\'e Saint-Jean, University of Alberta}\\
{\it Edmonton (AB), Canada}\\
{\it pjgaudre@ualberta.ca}\\
$^{\S}${\it Facult\'e Saint-Jean, University of Alberta}\\
{\it Edmonton (AB), Canada}\\
{\it hsafouhi@ualberta.ca}}

\date{}

\maketitle

% {\bf AMS classification:} \hskip 0.15cm 65L10, 65L20

\textbf{Abstract.}

The present contribution concerns the computation of energy eigenvalues of a perturbed anharmonic coulombic potential with irregular singularities using a combination of the Sinc collocation method and the double exponential transformation. This method provides a highly efficient and accurate algorithm to compute the energy eigenvalues of one-dimensional time-independent Schr\"odinger equation. The numerical results obtained illustrate clearly the highly efficiency and accuracy of the proposed method. All our codes are written in Julia and are available on github at \url{https://github.com/pjgaudre/DESincEig.jl}.

\vspace*{0.5cm}

Schr\"odinger equation, anharmonic potentials, inverse power potentials, Sinc collocation method, double exponential transformation.

\clearpage
\section{Introduction}

This paper continues the series of previous studies~\cite{Safouhi49, Safouhi52, Safouhi53} concerning the accurate and efficient computation of energy eigenvalues of one-dimensional time-independent Schr\"odinger equation with the anharmonic potential, using the combination of the Sinc collocation method (SCM)~\cite{Stenger-23-165-81, Stenger-121-379-00, Jarratt1990a, Alquran2010, Eggert-Jarratt-Lund-69-209-87} and the double exponential (DE) transformation~\cite{Takahasi-Mori-9-721-74, Sugihara2004, Mori-Sugihara-127-287-01, Sugihara2002}. We refer to this method as DESCM.

The Sinc collocation methods (SCM) consists of expanding the solution of Schr\"odinger equation using a basis of Sinc functions. By evaluating the resulting approximation at the Sinc collocation points, we obtain a generalized matrix eigenvalue problem for which the resulting eigenvalues are approximations to the energy eigenvalues of the Schr\"odinger equation.

The SCM have been used extensively during the last 30 years to solve many problems in numerical analysis~\cite{Stenger-23-165-81, Stenger-121-379-00} and it has been shown that it yields exponential convergence.

The DE transformation is a conformal mapping which allows for the function being approximated by a Sinc basis to decay double exponentially at both infinities. The effectiveness of the DE transformnation has been the subject of numerous studies~\cite{Sugihara2004, Mori-Sugihara-127-287-01, Sugihara2002} and it has been shown that DE yields the best available convergence properties for problems with end point singularities or infinite domains.

In~\cite{Safouhi49}, we applied DESCM to approximate the energy eigenvalues for anharmonic oscillators given by:
$$
V(x) = \sum_{i=1}^{m} c_{i} x^{2i} \quad \textrm{with} \quad c_{m}>0.
$$

The above expression generalizes potentials that have been previously treated~\cite{Bender-Orszag-78, Weniger-246-133-96, Weniger-Cizek-Vinette-34-571-93, Zamastil-Cizek-Skala-276-39-99, Patnaik-35-1234-87, Adhikari-Dutt-Varshni-131-217-88, Datta-Rampal-23-2875-81}, namely:
$$
\left\{
\begin{array}{rclll}
V(x) & = & x^{2}+\beta \, x^{2m} & \quad \textrm{for} \quad & \beta > 0 \quad \textrm{and} \quad m=2,3,4.\ldots\\[0.25cm]
V(x) & = & \lambda \, x^{2m}   & \quad \textrm{for} \quad & m=2,3,...\\[0.25cm]
V(x) & = & \mu \, x^{2}+\lambda \, x^{4}+\eta \, x^{6}  &  \quad \textrm{with} \quad & \eta >0\\[0.25cm]
V(x) & = & \frac{1}{2}\, x^{2}+\lambda \, x^{4}  & \quad \textrm{as} \quad & \lambda \rightarrow 0^{+} \qquad \textrm{and} \qquad \lambda \rightarrow \infty\\[0.25cm]
V(x) & = & a \, x^{2}+b\, x^{4}+c \, x^{6}   & \quad \textrm{as} \quad  & c \rightarrow \infty.
\end{array}\right.
$$

In~\cite{Safouhi52}, we introduced an algorithm based on DESCM  for the numerical treatment of the harmonic oscillators perturbed by rational functions in their general form given by:
$$
V(x) = \omega \, x^{2m} + \frac{p(x)}{q(x)}.
$$

Recently~\cite{Safouhi53}, we presented an algorithm based on DESCM for the coulombic potentials given by:
$$
V(x)  =  \frac{a_{-2}}{x^2}+\frac{a_{-1}}{x} + \sum_{i=1}^n a_i x^i.
$$

In~\cite{Safouhi49, Safouhi52, Safouhi53}, we have presented numerous advantages of DESCM over the existing alternatives. DESCM is not case specific and is insensitive to changes in the potential parameters. DESCM can be applied to a large set of anharmonic, rational and coulombic potentials. We have also demonstrated the high efficiency and accuracy of the method when dealing with multiple-well potentials. In addition, the method has a near-exponential convergence rate and the matrices generated by the DESCM have useful symmetric properties which simplify considerably the computation of their eigenvalues.

In the present contribution, we demonstrate that DESCM can also be applied to compute eigenvalues of the perturbed coulombic potential with irregular singulatrities given in their general form by:
$$
V(x) = \displaystyle \sum_{i=-p}^q a_i x^i  \qquad \textrm{with} \qquad p \geq 3 \quad \textrm{and} \quad a_{-p} >0.
$$

For the past three decades, the numerical resolution of the Schr\"odinger equation in the context of anharmonic oscillators has been approached mainly by case-specific methods, see for example \cite{Bender-Wu-184-1231-69, Bender-Orszag-78, Weniger-246-133-96, Weniger-Cizek-Vinette-34-571-93, Zamastil-Cizek-Skala-276-39-99, Patnaik-35-1234-87, Burrows1989, Adhikari-Dutt-Varshni-131-217-88, Datta-Rampal-23-2875-81,Tater-20-2483-87}.

In the case of potentials of the form:
$$
V(r) = ar^2 + br^{-2} + cr^{-4} + dr^{-6},
$$
an exact closed-form solution of the two-dimensional Schr\"odinger equation was obtained by applying an ansatz to the eigenfunction. Restrictions on the parameters of the given potential and the angular momentum were obtained \cite{Dong2000}. Another study was carried out in the case of a central potential:
$$
V(r) = \displaystyle{\sum_{i=1}^{4} a_{i} r^{-i}},
$$
where $V(r)$ is proportional to $r^{-\beta}$ in an arbitrary number of dimensions. The presence of a single term in the potential makes it impossible to use existing algorithms, which only work for quasi-exactly-solvable problems. Nevertheless, the analysis of the stationary Schr\"odinger equation in the neighborhood of the origin and at infinity provides relevant information about the desired solutions for all values of the radial coordinates. The original eigenvalue equation is mapped into a differential equation with milder singularities, and the role played by the particular case $\beta = 4$ is elucidated. In general, whenever the parameter $\beta$ is even and greater than $4$, a recursive algorithm for the evaluation of eigenfunctions is obtained. In the two-dimensional case, the exact form of the ground-state wave function is obtained for a potential containing a finite number of inverse powers of r, with the associated energy eigenvalue \cite{Dong1999}.

In \cite{Fernandez1991a} exact solutions for the quantum-mechanical harmonic oscillator were obtained when dealing with a perturbation potential belonging to the class of polynomial functions of $1/r$.
It has been shown that for some of the eigenfunctions, it is possible to calculate the expectation values in closed form. These eigenfunctions, therefore are suitable trial functions for the application of the variational method to related nonsolvable problems.

The eigenvalues of the potentials~\cite{Gonul2001}:
$$
V_1(r) =  \displaystyle{\sum_{i=1}^{4}a_{i} r^{-i}} \qquad \textrm{and} \qquad V_2(r)  = \displaystyle{ \sum_{i=-1}^{3} b_{i} r^{-2i}},
$$
were obtained in $N$-dimensional space for special cases, such as the Kratzer and Goldman-Krivchenkov potentials. It has been found that these potentials, which had not previously been linked, are explicitly dependent in higher-dimensional spaces~\cite{Gonul2001}. By using the power series method via a suitable ansatz to the wave function with parameters that may possibly exist in the potential function for the first time, an exact solution for the radial Schr\"odinger equation in N-dimensional Hilbert space was obtained for the potential $V_1(r)$~\cite{Khan2009}.

Exact analytical expressions for the energy spectra and potential parameters were obtained in terms of linear combinations of known parameters of radial quantum number $n$, angular momentum quantum number $l$, and the spatial dimension $N$. The coefficients of the expansion of the wave function ansatz were generated through the two-term recurence relation for odd/even solutions.

In \cite{Landtman1993}, a numerical method using a basis set consisting of $B$-splines was used to approximate the radial Schr\"odinger equation with the following anharmonic potential:
$$
V(r) = a \, r^2 + b \, r^{-4} + c \, r^{-6},
$$
and accurate approximations of the energy eigenvalues were obtained for the lowest few states. However, in the case of the first excited state, the values obtained deviate significantly from the reference values presented in~\cite{Kaushal-Parashar-1992}.

In \cite{Varshni1993}, it has been shown that the wavefunction proposed in~\cite{Kaushal-Parashar-1992} has four different possible sets of solutions for the wavefunction parameters and the constraint between the potential parameters. These four solutions were obtained analytically. The wavefunctions corresponding to these sets of parameters represent one the following states: the ground state, the first excited state, or the second excited state, depending on the values of the parameters of the wavefunction. It has also been shown that the values:
$$
a = 1, \quad b = - 11.25  \quad \textrm{and} \quad c = 3.515625,
$$
satisfy the constraint between the parameters both for the ground state and the first excited state.

In \cite{Lopez-Ortega2015}, two conditionally exactly solvable inverse power law potentials of the form:
$$
V(r) = \displaystyle{\sum_{i=1}^{n} a_{i} r^{-\alpha i}},
$$
whose linearly independent solutions include a sum of two confluent hypergeometric functions were given. Furthermore, it is noted that they are partner potentials and multiplicative shape invariant. The method used to find the solutions works with the two Schr\"odinger equations of the partner potentials.

In~\cite{Mikulski2015}, exact eigenvalues were found form the radial Schr\"odinger equation for the anharmonic potential given by:
$$
V(r) =\displaystyle{ \sum_{i=1}^{4} a_{i} r^{-i}}.
$$

In the study of the corresponding radial Schr\"odinger equation, all analytical calculations employ the mathematical formalism of supersymmetric quantum mechanics~\cite{Mikulski2015}. The novelty of this study is underlined by the fact that for the first time, recurrence formulas for rovibrational bound energy levels have been derived by employing factorization methods and an algebraic approach. Both the ground state and the excited states were determined by means of the hierarchy of the isospectral Hamiltonians. The Riccati nonlinear differential equation with superpotentials has been solved analytically. It has been shown that exact solutions exist when the potential and superpotential parameters satisfy certain supersymmetric constraints. The results obtained can be used both in computations of quantum chemistry and theoretical spectroscopy of diatomic molecules.

In a similar fashion to the use the ansatz for the eigenfunction, exact analytical solutions of the radial Schr\"{o}dinger equation have been obtained for two-dimensional pseudoharmonic and Kratzer potentials~\cite{Ozcelik1991a}. The bound-state solutions are easily calculated from the eigenfunction ansatz. The corresponding normalized wavefunctions are also obtained. By adopting the strategy of the ansatz for the eigenfunction  \cite{Simsek1994}, one can once more compute exact solutions of the Schr\"{o}dinger equation for two types of generalized potentials, namely the general Laurent type and four-parameter exponential potentials which can be reduced to the well-known types by choosing appropriate values for the parameters of:
$$
V(r)= \displaystyle{\sum_{i=-2j}^{2M} A_i r^{i}}   \quad \textrm{with} \quad  j, M \geq 0.
$$

From this overview of the literature on the numerical evaluation of the energy eigenvalues of Schr\"{o}dinger equation involving perturbed coulombic potentials with irregular singulatrities, it can be seen that generally, the proposed methods provide accurate results only for a specific class of potentials. The ideal would be to find a general method that is able to efficiently compute approximations of eigenvalues to a high pre-determined accuracy for more general potentials.

\clearpage
\section{General definitions and properties}\label{The double exponential Sinc collocation method}

The sinc function, analytic for all $z \in \mathbb{C}$ is defined by the following expression:
\begin{equation} \label{formula: sinc functions}
\textrm{sinc}(z) = \left\{ \begin{array}{cc} \dfrac{\sin(\pi z)}{\pi z} &\quad \textrm{for} \quad z \neq 0 \\[0.3cm]
1  &\quad \textrm{for} \quad z=0. \end{array} \right.
\end{equation}

The Sinc function $S(j,h)(x)$ for $h \in \mathbb{R}^{+}$ and $j \in \mathbb{Z}$ is defined by:
\begin{equation}\label{formula: Sinc function}
S(j,h)(x) = \sinc \left( \dfrac{x}{h}-j\right).
\end{equation}

The Sinc functions defined in \eqref{formula: Sinc function} form an interpolatory set of functions with the discrete orthogonality property:
\begin{equation}
S(j,h)(kh)  =  \delta_{j,k}   \qquad \textrm{for} \qquad j,k \in \mathbb{Z},
\end{equation}
where $\delta_{j,k}$ is the Kronecker delta function. In other words, at all the Sinc grid points $x_{k} = kh$, we have:
$$
C_{N}(v,h)(x) = v(x).
$$

It is possible to expand well-defined functions as series of Sinc functions. Such expansions are known as Sinc expansions or Whittaker Cardinal expansions.
\begin{definition}\cite{Stenger1981}
Given any function $v(x)$ defined everywhere on the real line and any $h>0$, the Sinc expansion of $v(x)$ is defined by the following series:
\begin{equation}\label{formula: Sinc expansion}
C(v,h)(x) = \sum_{j=-\infty}^{\infty} v_{j} S(j,h)(x),
\end{equation}
where $v_{j} = v(jh)$.
The symmetric truncated Sinc expansion of the function $v(x)$ is defined by the following series:
\begin{equation}\label{formula: Truncated Sinc expansion}
C_{N}(v,h)(x) = \sum_{j=-N}^{N} v_{j} \, S(j,h)(x) \qquad \textrm{for} \qquad N \in \mathbb{N}.
\end{equation}
\end{definition}

In~\cite{Stenger1981}, Stenger introduced a class of functions which are successfully approximated by a Sinc expansion. We present the definition for this class of functions bellow.
\begin{definition} \cite{Stenger1981} \label{defintion: Bd function space}
Let $d>0$ and let $\mathscr{D}_{d}$ denote the strip of width $2d$ about the real axis:
\begin{equation}
\mathscr{D}_{d} = \{ z \in \mathbb{C} : |\,\Im (z)|<d \}.
\end{equation}
For $\epsilon \in(0,1)$, let $\mathscr{D}_{d}(\epsilon)$ denote the rectangle in the complex plane:
\begin{equation}
\mathscr{D}_{d}(\epsilon) = \{z \in \mathbb{C} : |\,\Re(z)|<1/\epsilon, \, |\,\Im (z)|<d(1-\epsilon) \}.
\end{equation}
Let ${\bf B}_{2}(\mathscr{D}_{d})$ denote the family of all functions $g$ that are analytic in $\mathscr{D}_{d}$, such that:
\begin{equation}\label{formula: integral imaginary}
\displaystyle \int_{-d}^{d} | \,g(x+iy)| \textrm{d}y \to 0 \;\; \textrm{as} \;\;  x \to \pm \infty \quad  \textrm{and} \quad \displaystyle \lim_{\epsilon \to 0} \left(  \int_{\partial \mathscr{D}_{d}(\epsilon)}  |g(z)|^{2} |\textrm{d}z| \right)^{1/2} <\infty.
\end{equation}
\end{definition}

A function $\omega(x)$ is said to decay double exponentially at infinities if there exist positive constants $A, B,\gamma$ such that:
\begin{equation}
|\, \omega(x)| \leq A \exp( -B \exp(\gamma|\,x|)) \qquad \textrm{for} \qquad x\in \mathbb{R}.
\end{equation}

The double exponential transformation is a conformal mapping $\phi(x)$ which allows for the solution of~\eqref{formula: transformed Schrodinger equation} to have double exponential decay at both infinities.

\section{The double exponential Sinc collocation method}

The Schr\"{o}dinger equation with semi-infinite zero boundary conditions is given by:
\begin{align}\label{formula: Schrodinger equation sturm liouville}
\left(- \dfrac{d^2}{dx^2} + V(x) \right) \, \psi(x) &  =   E \, \psi(x) \qquad \textrm{for} \qquad  0<x< \infty
\nonumber\\
\psi(0) =  \psi(\infty) & = 0,
\end{align}
and where $V(x)$ stands for the perturbed anharmonic coulombic potential of the general form given by:
\begin{equation}
V(x) = \displaystyle \sum_{i=-p}^q a_i x^i \quad \textrm{with} \quad p \geq 3 \quad \textrm{and} \quad a_{-p} >0.
\end{equation}

The above potentials, are called singular potentials~\cite{Frank1971}, and they have an irregular singular point of rank $\frac{p-2}{2}$ at $x = 0$.

Eggert et al.~\cite{Eggert1987} demonstrate that applying an appropriate substitution to the boundary value problem~\eqref{formula: Schrodinger equation sturm liouville}, results in a symmetric discretized system when using Sinc expansion approximations. This change of variable is given by~\cite{Eggert1987}:
\begin{equation}
v(x) = \left(\sqrt{ (\phi^{-1})^{\prime} } \, \psi \right) \circ \phi(x) \qquad \Longrightarrow \qquad  \psi(x)  =  \dfrac{ v \circ \phi^{-1}(x)}{\sqrt{ (\phi^{-1}(x))^{\prime}}},
\label{formula: EggertSub}
\end{equation}
where $\phi^{-1}(x)$ a conformal map of a simply connected domain in the complex plane with boundary points $a\neq b$ such that:
\begin{equation}
\phi^{-1}(a)=-\infty \qquad \textrm{and}  \qquad \phi^{-1}(b)=\infty.
\end{equation}

Applying the substitution~\eqref{formula: EggertSub} to~\eqref{formula: Schrodinger equation sturm liouville}, we obtain:
\begin{align}\label{formula: transformed Schrodinger equation}
\hat{\mathcal{H}} \, v(x) & =  - v^{\prime \prime}(x) + \tilde{V}(x) \, v(x)
\nonumber\\ &  = \, E (\phi^{\prime}(x))^{2} \, v(x) \qquad \qquad \textrm{with}  \qquad \lim_{ |x| \to \infty} v(x) = 0,
\end{align}
where:
\begin{equation}
\tilde{V}(x)  =  - \sqrt{\phi^{\prime}(x)} \, \dfrac{{\rm d}}{{\rm d} x} \left( \dfrac{1}{\phi^{\prime}(x)} \dfrac{{\rm d}}{{\rm d} x} \left( \sqrt{\phi^{\prime}(x)}\right)  \right) + (\phi^{\prime}(x))^{2} \, V(\phi(x)).
\end{equation}

To implement the DESCM, we start by approximating the solution of~\eqref{formula: transformed Schrodinger equation} by a truncated Sinc expansion~\eqref{formula: Truncated Sinc expansion}.
Inserting~\eqref{formula: Truncated Sinc expansion} into~\eqref{formula: transformed Schrodinger equation}, we obtain the following system of equations:
\begin{eqnarray}
\widehat{\mathcal{H}}\, C_{N}(v,h)(x_{k}) & = & \sum_{j=-N}^{N} \left(-\frac{d^{2}}{dx^{2}} S(j,h)(x_{k}) + \widetilde{V}(x_{k}) \, S(j,h)(x_{k}) \right)\,  v(x_{j})
\nonumber\\
& = & {\mathcal E}\, \displaystyle{\sum_{j=-N}^{N}}\left(S(j,h)(x_{k})\, (\phi'(x_{k}))^{2}\right)\, v(x_{j}) \qquad \textrm{for}  \qquad   k=-N,\ldots,N,
\label{discret1}
\end{eqnarray}
where $x_{j}=jh$ for $j=-N,\ldots,N$ are the collocation points and by ${\mathcal E}$ the approximation of the eigenvalue $E$ in \eqref{formula: transformed Schrodinger equation}.

Equation (\ref{discret1}) can be written as follows:
\begin{eqnarray}
\widehat{\mathcal{H}}\, C_{N}(v,h)(x_{k}) & = & \sum_{j=-N}^{N} \left(-\frac{1}{h^{2}} \, \delta^{(2)}_{j,k} + \widetilde{V}(kh) \, \delta^{(0)}_{j,k}\right)\, v(jh)
\nonumber\\
& = & {\mathcal E} \, \displaystyle{\sum_{j=-N}^{N}} \left(\delta^{(0)}_{j,k}(\phi'(kh))^{2}\right)\, v(jh) \qquad \textrm{for}  \qquad   k=-N,...,N,
\label{discret2}
\end{eqnarray}
where:
$$
\delta^{(l)}_{j,k}=h^{l} \, \left(\frac{d}{dx}\right)^{l} \, S(j,h)(x)\Big|_{x=kh},
$$
and  more precisely:
\begin{equation}
\delta _{j,k}^{\left( 2\right) }=
\left\{
\begin{array}{lll}
-\dfrac{\pi ^{2}}{3}   & \quad \textrm{if}\quad & j=k \\[0.35cm]
\dfrac{\left( -2\right) \left( -1\right) ^{k-j}}{\left( k-j\right) ^{2}} & \quad \textrm{if} \quad \quad & j\neq k
\end{array}
\right.
\qquad \textrm{and} \qquad
\delta _{j,k}^{\left( 0\right) }=
\left\{
\begin{array}{lll}
1   & \quad \textrm{if} \quad & j=k \\[0.5cm]
0   & \quad \textrm{if} \quad & j\neq k
\end{array}
\right.
\end{equation}

The matrix form associated with (\ref{discret2}) is given by:
\begin{equation}\label{eigendisc}
\widehat{\mathcal{H}}\, \mathbf{C_{N}}(v,h) \,=\, \mathbf{H}\, \mathbf{v} \,=\,  {\mathcal E}\, \mathbf{D}^{2}\,\mathbf{v}
\quad \Longrightarrow \quad (\mathbf{H} -  {\mathcal E}\,  \mathbf{D}^{2})\,\mathbf{v} \,=\, 0,
\end{equation}
where:
\begin{eqnarray*}
\mathbf{v} & = & (v(-Nh),...,v(Nh))^{T}\\[0.25cm]
\mathbf{C_{N}}(v,h) & = & (C_{N}(v,h)(-Nh),...,C_{N}(v,h)(Nh))^{T}.
\end{eqnarray*}

The $(2N+1)\times (2N+1)$ matrices $\mathbf{H}$ and $\mathbf{D}^{2}$ are given by:
\begin{equation}\label{HD2}
\left\{
\begin{array}{lll}
\mathbf{H}_{jk}     & = & -\dfrac{1}{h^{2}}\, \delta^{(2)}_{j,k} + \widetilde{V}(kh)\,\delta^{(0)}_{j,k}\\[0.5cm]
\mathbf{D}_{jk}^{2} & = & (\phi'(kh))^{2} \,\delta^{(0)}_{j,k}
\end{array}
\right.\qquad  \textrm{with} \qquad  -N \leq j,k \leq N.
\end{equation}

\begin{theorem}\label{theorem: convergence of eigenvalues} \cite{Safouhi50}
Let $E$ and $v(x)$ be an eigenpair of the transformed Schr\"odinger equation~\eqref{formula: transformed Schrodinger equation}. Assume there exist positive constants $A,\beta_{L},\beta_{R},\gamma_{L}$ and $\gamma_{R}$ such that:
\begin{equation}\label{formula: xi growth condtion}
|\,v(x)| \leq  A
\left\{
\begin{array}{ccc}
 \exp( - \beta_{L} \exp(\gamma_{L} |x|) & \quad \textrm{for} \quad &  x \in (-\infty,0] \\[0.25cm]
 \exp( - \beta_{R} \exp(\gamma_{R} |x|))& \quad \textrm{for} \quad &  x \in [0, \infty ).
\end{array}
\right.
\end{equation}
Moreover, let $\gamma = \max \{ \gamma_{L}, \gamma_{R}\}$. If:
\begin{itemize}
\item $v \in {\bf B}_{2}(\mathscr{D}_{d})$ with $d \leq \dfrac{\pi}{2\gamma}$
\item there exists $\delta>0$ such that $\tilde{q}(x)\geq \delta^{-1}$
\item the mesh size $h$ is chosen optimally~by:
\begin{equation}
h = \dfrac{W(\pi d \gamma m / \beta)}{\gamma m},
\end{equation}
where $n$ and $\beta$ are given by:
\begin{equation}\label{formula: collocation points M}
\begin{cases}
m = M , \quad \beta = \beta_{L} & \quad \textrm{if} \quad \quad \gamma_{L} > \gamma_{R} \\
m = N , \quad \beta = \beta_{R} & \quad \textrm{if} \quad \quad \gamma_{R} > \gamma_{L} \\
m = M , \quad \beta = \beta_{L} & \quad \textrm{if} \quad \quad \gamma_{L} = \gamma_{R} \quad \textrm{and} \quad \beta_{L} \geq \beta_{R} \\
m = N , \quad \beta = \beta_{R} & \quad \textrm{if} \quad \quad \gamma_{L} = \gamma_{R} \quad \textrm{and} \quad \beta_{R} \geq \beta_{L},
\end{cases}
\end{equation}
\end{itemize}
then there is an eigenvalue ${\mathcal E}$ of the generalized eigenvalue problem (\ref{eigendisc}) satisfying:
\begin{equation}
|{\mathcal E}-E|   \,\leq\,\, K_{v,d} \,\sqrt{\delta E} \,\left(\dfrac{m^{5/2}}{\log(m)^{2}} \right)  \,
\exp \left(-  \dfrac{\pi d \gamma m}{\log(\pi d \gamma m/\beta)} \right) \qquad \textrm{as} \qquad m \to \infty,
\end{equation}
where $W$ is the Lambert-$W$ function and $K_{v,d}$ is a constant that depends on $v$ and $d$.
\end{theorem}

\section{The perturbed anharmonic coulombic potential}

To implement the DE transformation, we must choose a conformal map $\phi$ which would result in a solution of the transformed Schr\"{o}dinger equation~\eqref{formula: EggertSub} that decays double exponentially.

Since the perturbed anharmonic coulombic potential is analytic in $\mathbb{C}\backslash \{0\}$ and grows to infinity as $x\to\infty$ and as $x \to 0$, the wave function is also analytic in $\mathbb{C}\backslash \{0\}$  and normalizable over~$\mathbb{R}^{+}$.

To find the decay rate of the solution of the Schr\"odinger equation with the  perturbed anharmonic coulombic  potential, we use the WKB method. Substituting the ansatz $\psi(x) = e^{S(x)}$ into the Schr\"odinger equation and simplifying, we obtain:
\begin{equation}\label{formula: asym schrodinger}
S^{\prime \prime}(x) + \left( S^{\prime}(x) \right)^2 - \sum_{i=-p}^{q} a_{i}x^{i} + E =0.
\end{equation}

Let us treat the case as $x \to \infty$ first.
Since $q>1$, we have:
$$
S^{\prime \prime}(x) = o \left( S^{\prime}(x)^2 \right) \qquad \textrm{as} \qquad x \to \infty.
$$

Hence:
\begin{equation}
\left( S^{\prime}(x) \right)^2  \sim a_{q} x^{q} \qquad \textrm{as} \qquad x \to \infty.
\end{equation}

Using the initial condition $\displaystyle \lim_{x \to \infty} \psi(x) = 0$, we obtain:
\begin{equation}\label{formula: 1 order asym}
S(x) \sim - \left(\dfrac{ 2\sqrt{a_{q}}\, }{q+2}\right) x^{(q+2)/2}  \qquad \textrm{as} \qquad x \to \infty.
\end{equation}

Let us now treat the case as $x \to 0^{+}$.  Since $p>2$, we have:
$$
S^{\prime \prime}(x) = o \left( S^{\prime}(x)^2 \right)  \qquad \textrm{as} \qquad x \to 0^{+}.
$$

Hence:
\begin{equation}
\left( S^{\prime}(x) \right)^2  \sim a_{-p} x^{-p} \qquad \textrm{as} \qquad x \to 0^{+}.
\end{equation}

Using the initial condition $\displaystyle \lim_{x \to 0^{+}} \psi(x) = 0$, we obtain:
\begin{equation}\label{formula: 2 order asym}
S(x) \sim - \left(\dfrac{ 2\sqrt{a_{-p}}\, }{p-2}\right) x^{-(p-2)/2}  \qquad \textrm{as} \qquad x \to 0^{+}.
\end{equation}

From~\eqref{formula: 1 order asym} and \eqref{formula: 2 order asym}, it follows that $\psi(x)$ has the following decay rates:
\begin{equation}\label{formula: wave function asymptotic}
\psi(x) = \begin{cases}
 {\cal O} \left(\exp \left[ - \left(\dfrac{ 2\sqrt{a_{q}}\, }{q+2}\right) x^{(q+2)/2} \right] \right) & \qquad \textrm{as} \qquad x \to \infty \\[0.5cm]
{\cal O} \left(\exp \left[ - \left(\dfrac{ 2\sqrt{a_{-p}}\, }{p-2}\right) x^{-(p-2)/2} \right] \right) & \qquad \textrm{as}  \qquad x \to 0^{+}.
\end{cases}
\end{equation}

Away from both boundary points, the wave function $\psi(x)$ undergoes oscillatory behaviour. As can be seen from~\eqref{formula: wave function asymptotic}, the wave function $\psi(x)$ decays only single exponentially. However, by using the conformal mapping $\phi(x)=e^{x}$, we have:
\begin{align}\label{formula: wave function asymptotic 2}
|v(x)| & =  \left| \dfrac{\psi \circ \phi(x)}{\sqrt{\phi^{\prime}(x)}} \right| \nonumber\\
& =  \left| e^{-x/2}\psi(e^{x}) \right| \nonumber\\[0.25cm]
& =
\begin{cases}
 {\cal O} \left(\exp \left[ - \left(\dfrac{ 2\sqrt{a_{q}}\, }{q+2}\right) \exp\left\{ \left(\dfrac{q+2}{2} \right)x\right\} \right] \right) & \qquad \textrm{as} \qquad x \to \infty \\[0.5cm]
{\cal O} \left(\exp \left[ - \left(\dfrac{ 2\sqrt{a_{-p}}\, }{p-2}\right) \exp\left\{ -\left(\dfrac{p-2}{2} \right)x\right\}  \right] \right) & \qquad \textrm{as}  \qquad x \to -\infty.
\end{cases}
\end{align}

According to Theorem~\ref{theorem: convergence of eigenvalues}, we have the following  parameters:
\begin{align}
\beta_{L} & = \dfrac{ 2\, \sqrt{a_{-p}}\, }{p-2} \qquad \textrm{and} \qquad \beta_{R}  = \dfrac{ 2\, \sqrt{a_{q}}\, }{q+2},
\end{align}
as well as
\begin{align}
\gamma_{L} & = \dfrac{p-2}{2} \qquad \textrm{and} \qquad
\gamma_{R}  = \dfrac{q+2}{2}.
\end{align}

Hence, the optimal mesh size $h$ is chosen by:
\begin{equation}
h = \dfrac{W(\pi d \gamma n / \beta)}{\gamma \, n},
\end{equation}
where $n$ and $\beta$ are given by:
\begin{equation}
\begin{cases}
n = M , \quad \beta = \beta_{L} & \quad \textrm{if} \quad \quad p-q > 4 \\
n = N , \quad \beta = \beta_{R} & \quad \textrm{if} \quad \quad p-q < 4 \\
n = M , \quad \beta = \beta_{L} & \quad \textrm{if} \quad \quad p-q = 4 \quad \textrm{and} \quad a_{-p} \geq a_{q} \\
n = N , \quad \beta = \beta_{R} & \quad \textrm{if} \quad \quad p-q= 4  \quad \textrm{and} \quad a_{-p} \leq a_{q}.
\end{cases}
\end{equation}
and $\gamma$ is given by:
\begin{equation}
\gamma = \max \left\{ \dfrac{p-2}{2}, \dfrac{q+2}{2} \right\}.
\end{equation}

To obtain an approximation of the eigenfunctions, we invert the variable transformation:
$$
v(x) =   e^{-x/2}\, \psi(e^{x}),
$$
to obtain:
\begin{eqnarray}
\psi(x) & = & \sqrt{x} \, v(\ln(x))
\nonumber \\ &\approx  & \sqrt{x} \sum_{j=-M}^{N} v_{j}\, S(j,h)(\ln(x)).
\label{formula: approximate wave function}
\end{eqnarray}

\section{Numerical discussion}
In this section, we present numerical results for the energy values of the anharmonic coulombic potentials. All calculations are performed using the programming language Julia~\cite{Bezanson2012} in double precision. A double-precision floating-point format is a computer number format that occupies 8 bytes (64 bits) in computer memory. On average, this corresponds to about $15$-$17$ significant decimal digits.

The eigenvalue solvers in Julia utilize the linear algebra package {\it LAPACK}~\cite{Anderson1999}. The matrices ${\bf H}$ and ${\bf D}^{2}$ are constructed using formulas (\ref{HD2}). The approxiamte wavefunctions are obtained using equation \eqref{formula: approximate wave function}. To produce our figures, we used the Julia package {\it PyPlot}. The saturation effect observed in all figures is merely a consequence of rounding errors resulting from this computer number format.

Exact solutions are presented in~\cite{Fernandez1991a} for potentials of the form:
\begin{equation}\label{formula: exact potential}
V(x) = x^2 +\left[\dfrac{k(k+1)}{2}\right]\left( \dfrac{2k-1}{x^4} + \dfrac{k(k+1)}{2x^6}  \right) \qquad \Longrightarrow \qquad E_{k} =2k+3, \quad k\geq 0.
\end{equation}
Explicitly, we define the absolute error as:
\begin{equation}
{\rm Absolute \, \, error} = \left| {\mathcal E}_{n}(m) - E_{n} \right| \qquad \textrm{for} \quad \left\{\begin{array}{ccc} m &= &1,2,3,\ldots \\
n & = & 0,1,2,\ldots. \end{array} \right.
\end{equation}

When the eigenvalues are not known exactly, we use an approximation to the absolute error given by:
\begin{equation}\label{formula: approximate absolute error}
\epsilon_{n}(m) = \left | {\mathcal E}_{n}(m)-{\mathcal E}_{n}(m-1) \right | \qquad \textrm{for} \quad \left\{\begin{array}{ccc} m &= &1,2,3,\ldots \\
n & = & 0,1,2,\ldots. \end{array} \right.
\end{equation}

\begin{figure}[!ht]
\begin{center}
\begin{tabular}{cccc} \includegraphics[width=1.1\textwidth]{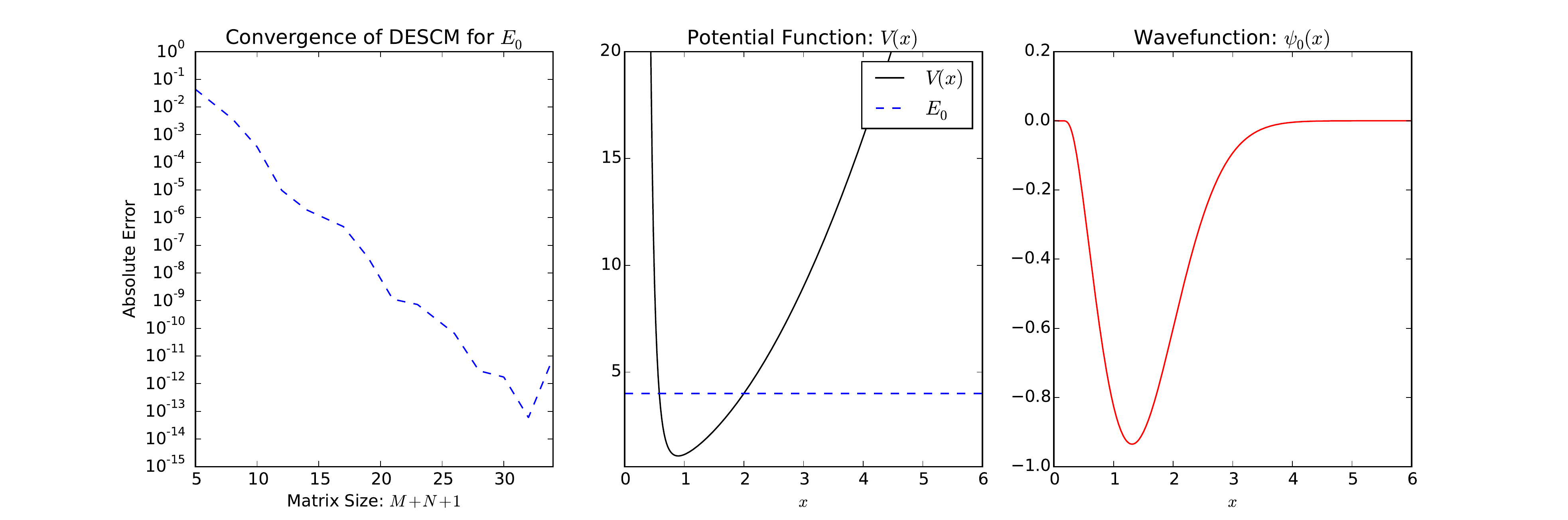}\\
(a) \\
\includegraphics[width=1.1\textwidth]{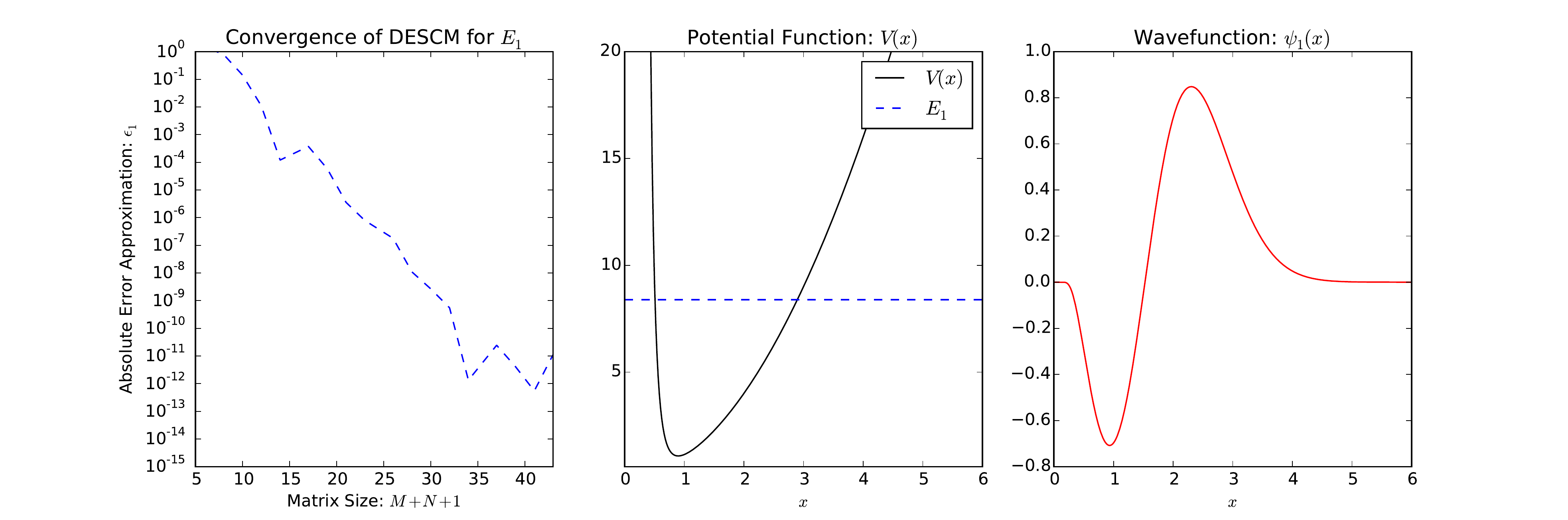}\\
(b)
\end{tabular}
\caption{Figure (a) displays the absolute error of the DESCM for the potential function $V(x)$ \eqref{formula: exact potential} for $k=1/2$ with exact eigenvalue $E_{0} = 4$, in addition to the associated wavefunction $\psi_{0}(x)$. Figure (b) displays the convergence of the DESCM for the same potential function $V(x)$ \eqref{formula: exact potential} for the first excited state in addition to the associated wavefunction $\psi_{1}(x)$.}
\label{figure: Exact Potential 0 and 1}
\end{center}
\end{figure}

In figure \ref{figure: Exact Potential 0 and 1}, we plotted the convergence of the DESCM for the ground and first excited states, the potential function in \eqref{formula: exact potential} with $k=1/2$, as well as the associated wavefunctions $\psi_{0}(x)$ and $\psi_{1}(x)$. As can be seen from the figure, the DESCM converges quickly in both cases, achieving over ten digits of accuracy for matrix sizes less than 35 by 35.

To demonstrate the robustness and power of the proposed method, we will now tackle more complicated potentials of the form:
\begin{equation}\label{Potential}
V(x) = \sum_{i=-p}^{q} a_{i}x^{i}.
\end{equation}

\begin{figure}[!ht]
\begin{center}
\begin{tabular}{cccc} \includegraphics[width=1.0\textwidth]{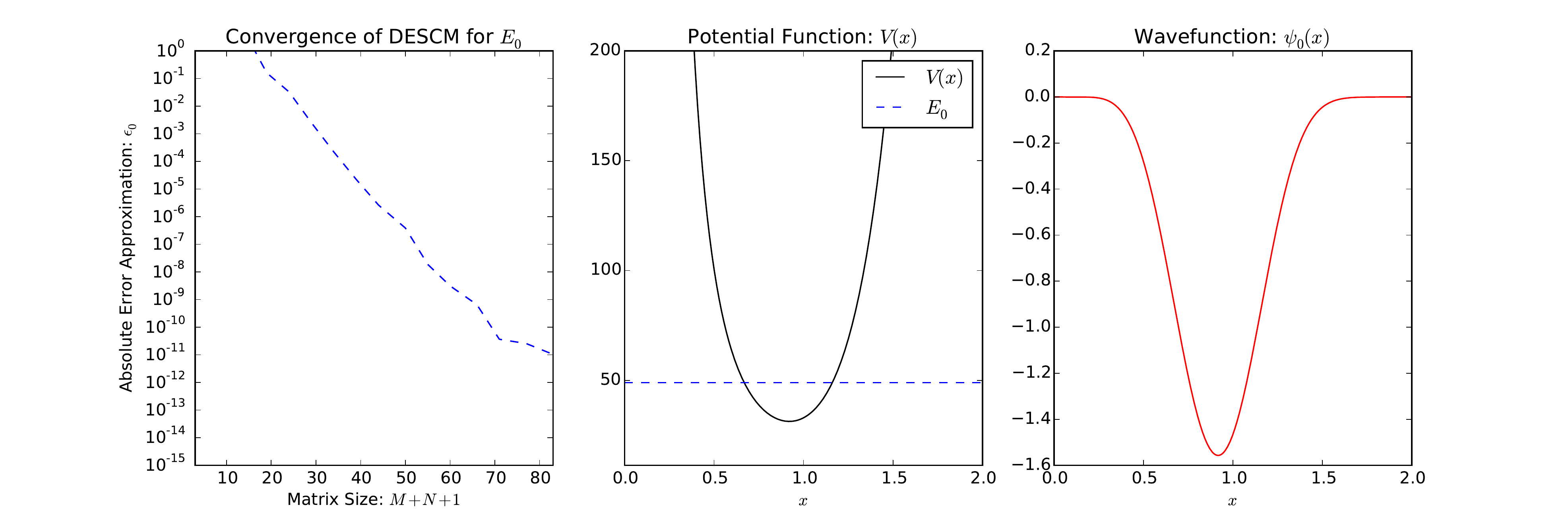}\\
(a) \\
\includegraphics[width=1.0\textwidth]{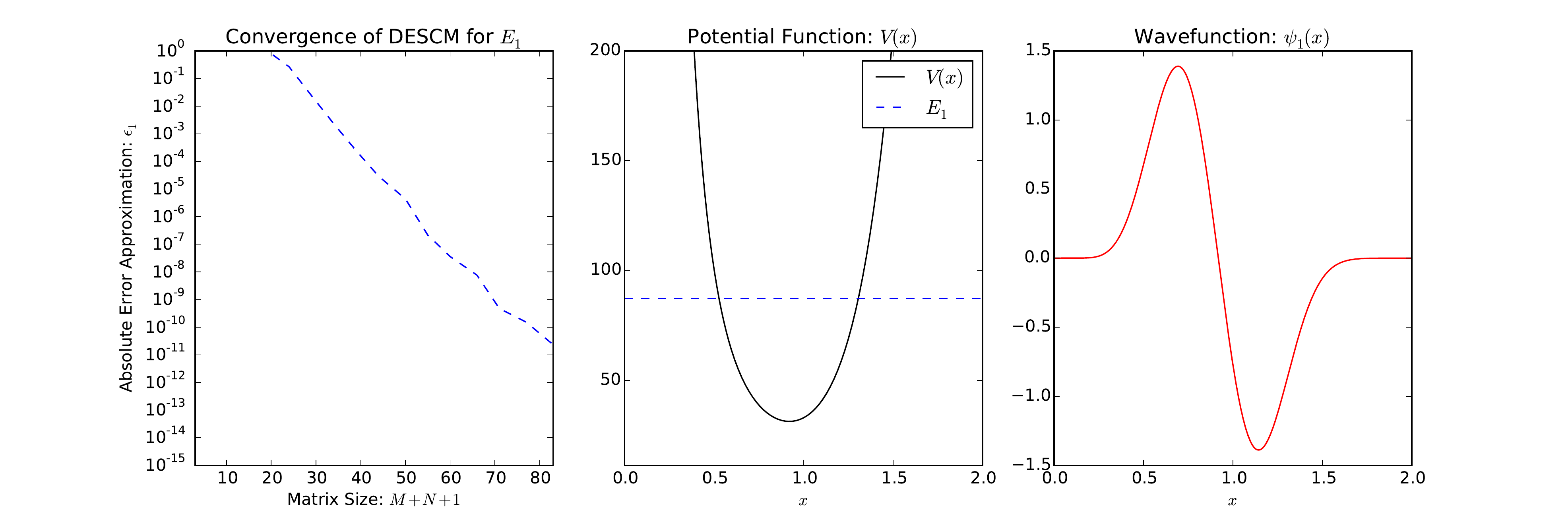}\\
(b) \\
\includegraphics[width=1.0\textwidth]{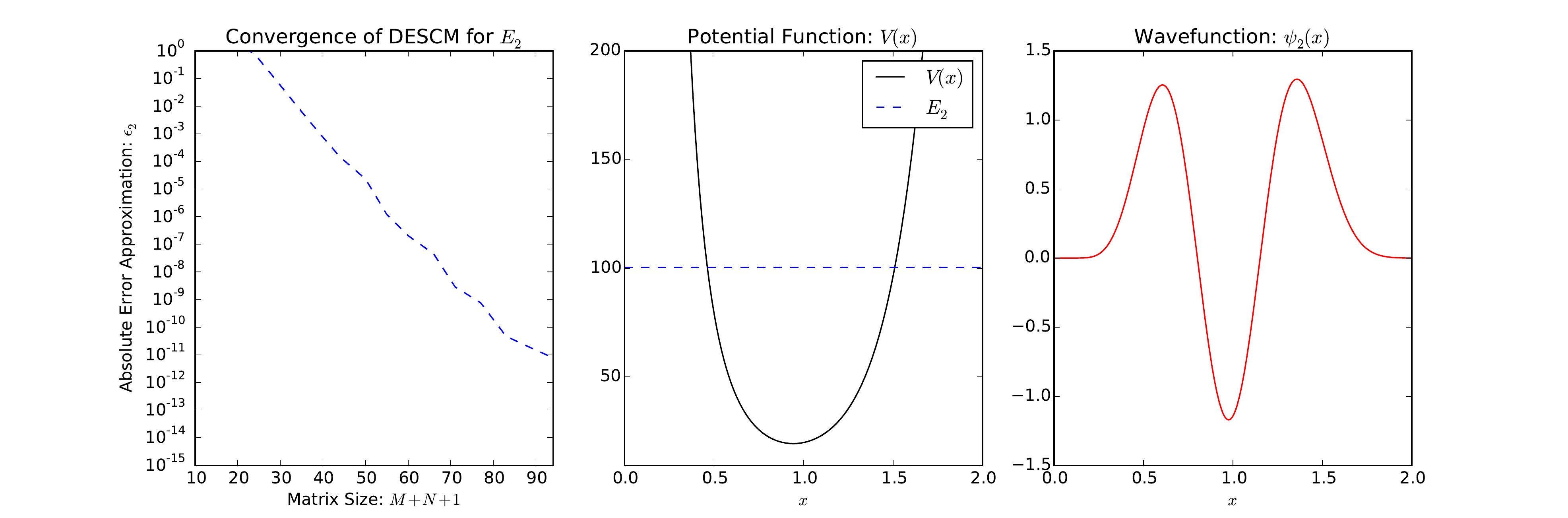}\\
(c)
\end{tabular}
\caption{Figure (a) displays the absolute error approximation of the DESCM for the potential function $V(x)=\sum_{i=-3}^{8} a_{i}x^{i}$, in addition to the associated wavefunction $\psi_{0}(x)$. Figure (b) displays the convergence of the DESCM for the same potential function $V(x)$ for the first excited state in addition to the associated wavefunction $\psi_{1}(x)$. Figure (c) displays the convergence of the DESCM for the same potential function $V(x)$ for the third excited state in addition to the associated wavefunction $\psi_{3}(x)$.}
\label{figure: Medium Potential}
\end{center}
\end{figure}

In figure \ref{figure: Medium Potential}, we present the convergence of the DESCM for the potential \eqref{Potential} with $p=3$ and $q=8$.
The $12$ coefficients $a_{i}$ involved in the expression of $V(x)$ are generated randomly using random number generator such that $-5\leq a_{i}\leq 5$ for $i= -2 ,\ldots,7$ . To assure positive values for the last two remaining values, we set $a_{-3} = 10.0$ and $ a_{8} = 5.0$.

\begin{figure}[!ht]
\begin{center}
\begin{tabular}{cccc} \includegraphics[width=1.05\textwidth]{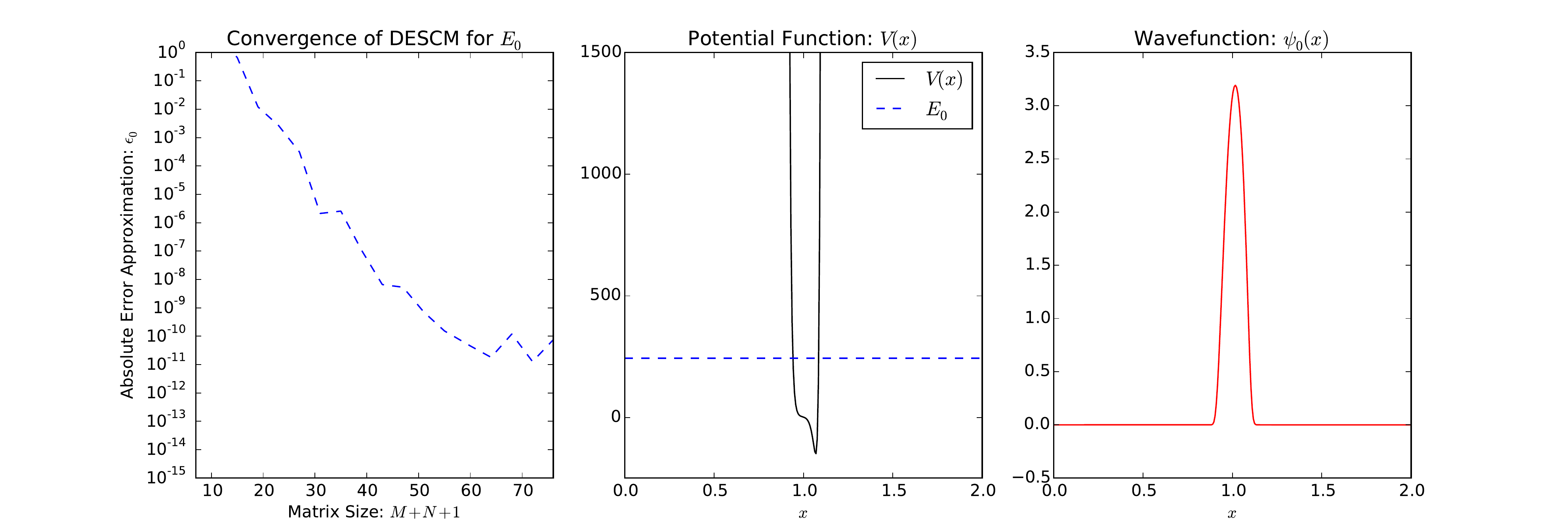}\\
(a) \\
\includegraphics[width=1.05\textwidth]{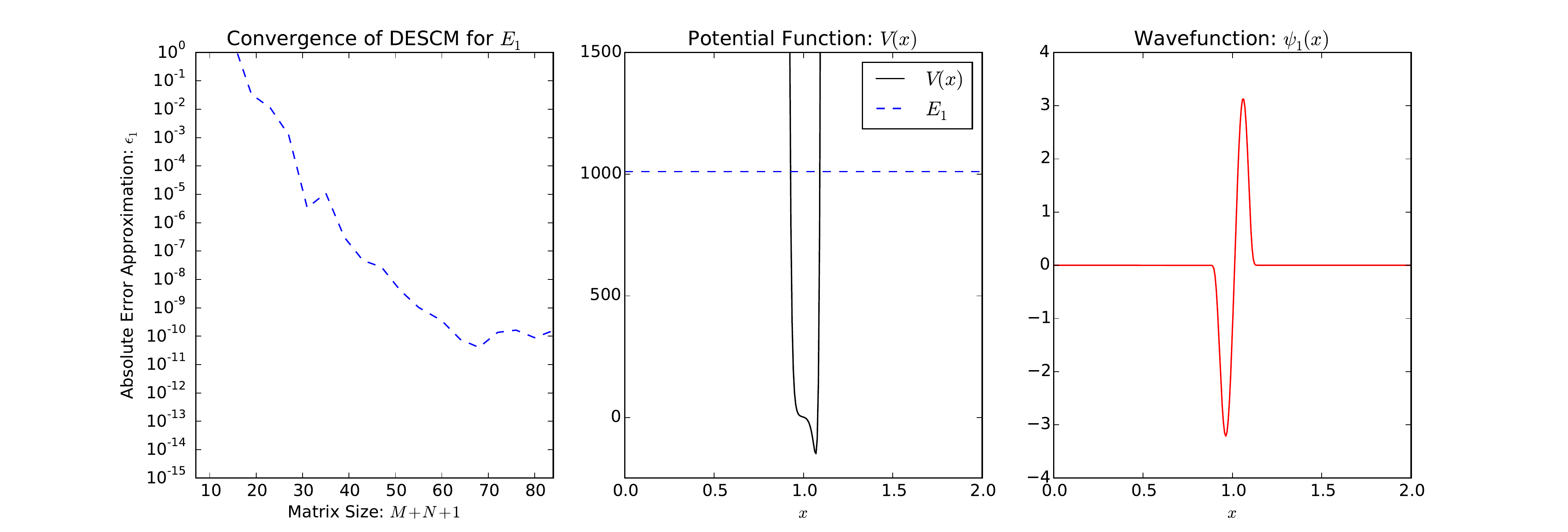}\\
(b) \\
\includegraphics[width=1.05\textwidth]{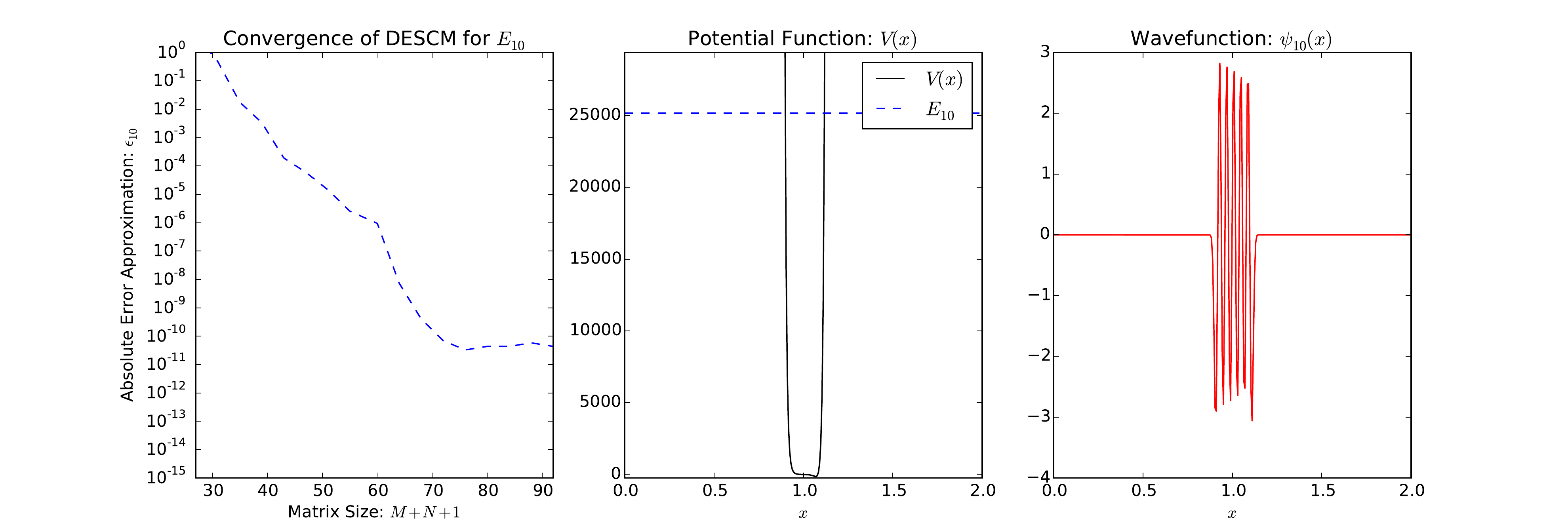}\\
(c)
\end{tabular}
\caption{Figure (a) displays the absolute error approximation of the DESCM for the potential function $V(x)=\sum_{i=-100}^{100} a_{i}x^{i}$, in addition to the associated wavefunction $\psi_{0}(x)$. Figure (b) displays the convergence of the DESCM for the same potential function $V(x)$ for the first excited state in addition to the associated wavefunction $\psi_{1}(x)$. Figure (c) displays the convergence of the DESCM for the same potential function $V(x)$ for the 10th excited state in addition to the associated wavefunction $\psi_{10}(x)$.}
\label{figure: Extreme Potential}
\end{center}
\end{figure}

In figure \ref{figure: Extreme Potential}, we present the convergence of the DESCM for the potential \eqref{Potential} with $p=q=100$.
Here again, we used a random number generator for the 201 coefficients $a_{i}$ such that $-1\leq a_{i}\leq 1$ for $i= -99 ,\ldots,99$. To assure positive values for the last two remaining values, we set $a_{-100} = a_{100} = 1.0$.

\section{Conclusion}
In the present contribution, we applied the DESCM to the Schr\"odinger equation with an anharmonic coulombic potential, which have the form of a Laurent series. The DESCM approximates the wave function of a transformed Schr\"odinger equation \eqref{formula: transformed Schrodinger equation} by as a Sinc expansion. By summing over $N+M+1$ collocation points, the implementation of the DESCM leads to a generalized eigenvalue problem with symmetric, positive definite matrices. In addition, we also state that the convergence of the DESCM in this case can be improved to the rate ${\cal O} \left( \left(\dfrac{n^{5/2}}{\log(n)^{2}} \right)  \exp \left(- \kappa^{\prime}\dfrac{n}{\log(n)} \right)\right)$ as $n\to \infty$, where $n$ is related to the dimension of the resulting generalized eigenvalue system and $\kappa^\prime$ is a constant that depends on the potential. As demonstrated in the numerical section, our application of this method on extreme potentials with a large number of coefficients was highly successful.

\section{Tables and Figures}

\clearpage
%\bibliography{/Users/hsafouhi/GoogleDrive/AMybibliography/My_Bibliography,/Users/hsafouhi/GoogleDrive/AMybibliography/library}

\begin{thebibliography}{10}

\bibitem{Safouhi49}
P.~Gaudreau, R.M. Slevinsky, and H.~Safouhi.
\newblock The double exponential sinc-collocation method for computing energy
  levels of anharmonic oscillators.
\newblock {\em Annals of Physics}, 360:520--538, 2015.

\bibitem{Safouhi52}
P.~Gaudreau and H.~Safouhi.
\newblock Double exponential sinc-collocation method for solving the energy
  eigenvalues of harmonic oscillators perturbed by a rational function.
\newblock {\em Journal of Mathematical Physics}, 58:101509 (1--15), 2017.

\bibitem{Safouhi53}
P.~Gaudreau T.~Cassidy and H.~Safouhi.
\newblock On the computation of eigenvalues of the anharmonic coulombic
  potential.
\newblock {\em Journal Mathematical Chemistry}, 56:477--492, 2017.

\bibitem{Stenger-23-165-81}
F.~Stenger.
\newblock Numerical methods based on {Whittaker} cardinal, or {Sinc} functions.
\newblock {\em SIAM Rev.}, 23:165--224, 1981.

\bibitem{Stenger-121-379-00}
F.~Stenger.
\newblock Summary of {Sinc} numerical methods.
\newblock {\em Journal of Computational and Applied Mathematics}, 121:379--420,
  2000.

\bibitem{Jarratt1990a}
M.~Jarratt, J.~Lund, and K.L. Bowers.
\newblock {Galerkin schemes and the Sinc-Galerkin method for singular
  Sturm-Liouville problems}.
\newblock {\em Journal of Computational Physics}, 89(1):41--62, 1990.

\bibitem{Alquran2010}
M.T. Alquran and K.~Al-Khaled.
\newblock {Approximations of Sturm-Liouville eigenvalues using Sinc-Galerkin
  and differential transform methods}.
\newblock {\em Applications and Applied Mathematics: An International Journal},
  5(1):128--147, 2010.

\bibitem{Eggert-Jarratt-Lund-69-209-87}
N.~Eggert, M.~Jarratt, and J.~Lund.
\newblock Sinc function computation of the eigenvalues of {Sturm-Liouville}
  problems.
\newblock {\em Journal of Computational Physics}, 69:209--229, 1987.

\bibitem{Takahasi-Mori-9-721-74}
H.~Takahasi and M.~Mori.
\newblock Double exponential formulas for numerical integration.
\newblock {\em RIMS}, 9:721--741, 1974.

\bibitem{Sugihara2004}
M.~Sugihara and T.~Matsuo.
\newblock {Recent developments of the Sinc numerical methods}.
\newblock {\em Journal of Computational and Applied Mathematics},
  164-165(1):673--689, 2004.

\bibitem{Mori-Sugihara-127-287-01}
M.~Mori and M.~Sugihara.
\newblock The double-exponential transformation in numerical analysis.
\newblock {\em Journal of Computational and Applied Mathematics}, 127:287--296,
  2001.

\bibitem{Sugihara2002}
M.~Sugihara.
\newblock {Double exponential transformation in the Sinc-collocation method for
  two-point boundary value problems}.
\newblock {\em Journal of Computational and Applied Mathematics},
  149(1):239--250, 2002.

\bibitem{Bender-Orszag-78}
C.M. Bender and S.A. Orszag.
\newblock {\em Advanced mathematical methods for scientists and engineers}.
\newblock Springer-Verlag New York, New York, 1978.

\bibitem{Weniger-246-133-96}
E.J. Weniger.
\newblock A convergent renormalized strong coupling perturbation expansion for
  the ground state energy of the quartic, sextic, and octic anharmonic
  oscillator.
\newblock {\em Ann. Phys. (NY)}, 246:133--165, 1996.

\bibitem{Weniger-Cizek-Vinette-34-571-93}
E.J. Weniger, J.~C\'{i}zek, and F.~Vinette.
\newblock The summation of the ordinary and renormalized perturbation series
  for the ground state energy of the quartic, sextic, and octic anharmonic
  oscillators using nonlinear sequence transformations.
\newblock {\em J. Math. Phys.}, 34:571--609, 1993.

\bibitem{Zamastil-Cizek-Skala-276-39-99}
J.~Zamastil, J.~C\'{i}zek, and L.~Sk\'{a}la.
\newblock Renormalized perturbation theory for quartic anharmonic oscillator.
\newblock {\em Ann. Phys. (NY)}, 276:39--63, 1999.

\bibitem{Patnaik-35-1234-87}
P.K. Patnaik.
\newblock Rayleigh-{Schr\"{o}dinger} perturbation theory for the anharmonic
  oscillator.
\newblock {\em Physical Review D}, 35:1234--1238, 1987.

\bibitem{Adhikari-Dutt-Varshni-131-217-88}
R.~Adhikari, R.~Dutt, and Y.P. Varshni.
\newblock On the averaging of energy eigenvalues in the supersymmetric {WKB}
  method.
\newblock {\em Physics Letters A}, 131:217--221, 1988.

\bibitem{Datta-Rampal-23-2875-81}
K.~Datta and A.~Rampal.
\newblock Asymptotic series for wave functions and energy levels of doubly
  anharmonic oscillators.
\newblock {\em Physical Review D}, 23:2875--2883, 1981.

\bibitem{Bender-Wu-184-1231-69}
C.M. Bender and T.T. Wu.
\newblock Anharmonic oscillator.
\newblock {\em The Physical Review}, 184:1231--1260, 1969.

\bibitem{Burrows1989}
B.L. Burrows, M.~Cohen, and T.~Feldmann.
\newblock {A unified treatment of Schrodinger's equation for anharmonic and
  double well potentials}.
\newblock {\em Journal of Physics A: Mathematical and General},
  22(9):1303--1313, 1989.

\bibitem{Tater-20-2483-87}
M.~Tater.
\newblock The {Hill} determinant method in application to the sextic
  oscillator: limitations and improvement.
\newblock {\em J. Phys. A: Math. Gen.}, 20:2483--2495, 1987.

\bibitem{Fernandez2008}
F.M. Fern{\'{a}}ndez.
\newblock {Calculation of bound states and resonances in perturbed Coulomb
  models}.
\newblock {\em Physics Letters A}, 372(17):2956--2958, 2008.

\bibitem{Fernandez1991}
F.M. Fern{\'{a}}ndez.
\newblock {Convergent power-series solutions to the Schr{\"o}dinger equation
  with the potential}.
\newblock {\em Physics Letters A}, 160(2):116--118, 1991.

\bibitem{Fernandez-Ma-Tipping-40-6149-89}
F.M. Fernandez, Q.~Ma, and R.H. Tipping.
\newblock Eigenvalues of the {Schr\"{o}dinger} equation via the
  {Riccati-Pad\'e} method.
\newblock {\em Physical Review A}, 40:6149--6153, 1989.

\bibitem{Dong2000}
S.~Dong.
\newblock {Exact solutions of the two-dimensional Schrodinger equation with
  certain central potentials}.
\newblock {\em International Journal of Theoretical Physics}, 39(4):1119--1128,
  2000.

\bibitem{Dong1999}
S.~Dong, Z.~Ma, and G.~Esposito.
\newblock {Exact solutions of the Schr{\"{o}}dinger equation with inverse-power
  potential}.
\newblock {\em Foundations of Physics Letters}, 12(5):11, 1999.

\bibitem{Fernandez1991a}
F.M. Fern{\'{a}}ndez.
\newblock {Exact and approximate solutions to the Schr{\"{o}}dinger equation
  for the harmonic oscillator with a singular perturbation}.
\newblock {\em Physics Letters A}, 160(6):511--514, 1991.

\bibitem{Gonul2001}
B.~Gonul, O.~Ozer, M.~Kocak, D.~Tutcu, and Y.~Cancelik.
\newblock {Supersymmetry and the relationship between a class of singular
  potentials in arbitrary dimensions}.
\newblock {\em Journal of Physics A: Mathematical and General}, 34:8271--8279,
  2001.

\bibitem{Khan2009}
G.R. Khan.
\newblock {Exact solution of N-dimensional radial Schr{\"{o}}dinger equation
  for the fourth-order inverse-power potential}.
\newblock {\em The European Physical Journal D}, 53(2):123--125, jun 2009.

\bibitem{Landtman1993}
M.~Landtman.
\newblock {Calculation of low lying states in the potential $V(r) = a r^2 + b
  r^{-4} + c r^{-6}$ using B-spline basis sets}.
\newblock {\em Physics Letters A}, 175(3-4):147--149, 1993.

\bibitem{Kaushal-Parashar-1992}
R.S. Kaushal and D~Parashar.
\newblock On the quantum bound states for the potential
  $v(r)=ar^2+br^{−4}+cr^{−6}$ using b-spline basis sets.
\newblock {\em Physics Letters A}, 170(5):335--338, 1992.

\bibitem{Varshni1993}
Y.P. Varshni.
\newblock {The first three bound states for the potential $V(r) = a r^2 + b
  r^{-4} + c r^{-6}$}.
\newblock {\em Physics Letters A}, 183(1):9--13, nov 1993.

\bibitem{Lopez-Ortega2015}
A.~L{\'{o}}pez-Ortega.
\newblock {New conditionally exactly solvable inverse power law potentials}.
\newblock {\em Physica Scripta}, 90(8):085202, aug 2015.

\bibitem{Mikulski2015}
Damian Mikulski, Jerzy Konarski, Krzysztof Eder, Marcin Molski, and Stanislaw
  Kabacinski.
\newblock {Exact solution of the Schr{\"{o}}dinger equation with a new
  expansion of anharmonic potential with the use of the supersymmetric quantum
  mechanics and factorization method}.
\newblock {\em Journal of Mathematical Chemistry}, 53(9):2018--2027, 2015.

\bibitem{Ozcelik1991a}
S.~\"Ozcelik and M.~Simsek.
\newblock {Exact solutions of the radial Schr{\"{o}}dinger equation for
  inverse-power potentials}.
\newblock {\em Physics Letters A}, 152(3-4):145--150, 1991.

\bibitem{Simsek1994}
M.~Simsek and S.~\"Ozcelik.
\newblock {Bound state solutions of the Schr{\"{o}}dinger equation for
  reducible potentials: general Laurent series and four-parameter
  exponential-type potentials}.
\newblock {\em Physics Letters A}, 186(1-2):35--40, mar 1994.

\bibitem{Stenger1981}
F.~Stenger.
\newblock {Numerical Methods Based on Whittaker Cardinal, or Sinc Functions}.
\newblock {\em SIAM Review}, 23(2):165--224, 1981.

\bibitem{Frank1971}
W.M. Frank, D.J. Land, and R.M. Spector.
\newblock {Singular Potentials}.
\newblock {\em Reviews of Modern Physics}, 43(1):36--98, 1971.

\bibitem{Eggert1987}
N.~Eggert, M.~Jarratt, and J.~Lund.
\newblock {Sinc function computation of the eigenvalues of Sturm-Liouville
  problems}.
\newblock {\em Journal of Computational Physics}, 69(1):209--229, 1987.

\bibitem{Safouhi50}
P.~Gaudreau, R.M. Slevinsky, and H.~Safouhi.
\newblock The double exponential sinc collocation method for singular
  sturm-liouville problems.
\newblock {\em Journal of Mathematical Physics}, 57:043505 (1--19), 2016.

\bibitem{Bezanson2012}
J.~Bezanson, S.~Karpinski, V.B. Shah, and A.~Edelman.
\newblock {Julia: A Fast Dynamic Language for Technical Computing}.
\newblock arXiv(1209.5145):1--27, 2012.

\bibitem{Anderson1999}
E.~Anderson, Z.~Bai, C.~Bischof, S.~Blackford, J.~Demmel, J.~Dongarra, J.~{Du
  Croz}, A.~Greenbaum, S.~Hammarling, A.~McKenney, and D.~Sorensen.
\newblock {\em {{LAPACK} Users' Guide}}.
\newblock Society for Industrial and Applied Mathematics, Philadelphia, PA,
  third edition, 1999.

\end{thebibliography}

\end{document}